\documentclass[10pt]{article}

\usepackage{amssymb,latexsym,amsfonts,amsmath,amsthm,verbatim,eucal,hyperref}

\theoremstyle{plain}
\newtheorem*{thm*}{Theorem}
\newtheorem{lem}{Lemma}
\theoremstyle{remark}
\newtheorem*{rmk*}{Remark}
\newtheorem*{ack*}{Acknowledgement}
\def\P{\mathbb{P}}
\def\E{\mathbb{E}}
\def\R{\mathbb{R}}
\def\S{S^{n-1}}
\def\eps{\varepsilon}
\def\sign{\operatorname{sign}}
\def\Ave{\operatorname{Ave}}

\begin{document}

\title{Kahane-Khinchin type Averages}
\author{Omer Friedland\textsuperscript{1}}
\maketitle

\begin{abstract}
\noindent We prove a Kahane-Khinchin type result with a few random
vectors, which are distributed independently with respect to an
arbitrary log-concave probability measure on $\R^n$. This is an
application of small ball estimate and Chernoff's method, that has
been recently used in the context of Asymptotic Geometric Analysis
in \cite{afm1}, \cite{afm2}.
\end{abstract}

\section{Introduction}
\footnotetext[1]{omerfrie@post.tau.ac.il; address: School of
Mathematical Sciences, Tel Aviv University, Ramat Aviv, Tel Aviv
69978, Israel.}

The classical \textbf{Kahane's inequality} (cf. Kahane \cite{ka})
states that for any $1\le p<\infty$ there exists a constant $K_p>0$
such that
\begin{equation}\label{kahane}
\int_0^1\|\sum_{i=1}^n r_i x_i\|dt \le (\int_0^1\|\sum_{i=1}^n r_i
x_i\|^p dt)^{1/p} \le K_p \int_0^1\|\sum_{i=1}^n r_i x_i\|dt ~,
\end{equation}

holds true for every $n$ and arbitrary choice of vectors $x_1 ,
\ldots , x_n \in X$, where $X$ is a normed space with the norm $\|
\cdot \|$, and $r_i$ is a Rademacher function that is given by
$r_i(t)=\sign\sin(2^{i-1}\pi t)$, $i\ge1$ (for more information
about Rademacher functions, see, e.g. Milman and Schechtman
\cite{ms}, section 5.5). Kwapie\'n \cite{kw} proved that $K_p\sim
\sqrt{p}$ is a constant depending only on $p$.

Note that, in view of the definition of $r_i(t)$, these integrals in
(\ref{kahane}) can also be represented as averages (for $p\ge1$)
\begin{equation}\label{average}
\int_0^1\|\sum_{i=1}^n r_i(t)x_i\|^p dt = \Ave_{\pm}\|\sum_{i=1}^n
\pm x_i\|^p ~.
\end{equation}

Bourgain, Lindenstrauss and Milman \cite{blm} proved that, if $v_1,
\ldots ,v_n$ are unit vectors in a normed space $(\R^n,\| \cdot \|)$
then, for any $\eps>0$ there exists a constant $C(\eps)>0$, such
that $N=C(\eps)n$ random sign vectors $\eps(1), \ldots , \eps(N) \in
\{-1,1\}^n$ satisfy with probability greater than $1-e^{-cn}$ that
\begin{equation}\label{blm}
(1-\eps) |||x||| \le {1\over N} \sum_{j=1}^N\|\sum_{i=1}^n \eps(j)_i
x_i v_i\| \le (1+\eps) |||x||| ~,
\end{equation}
for every $x\in\R^n$, where $ |||x||| =\Ave_{\pm}\|\sum_{i=1}^n \pm
x_i v_i\|$.

Obviously, the norm $||| \cdot |||$ is an unconditional norm, i.e.,
it is invariant under the change of signs of the coordinates. As
(\ref{blm}) shows, it is sufficient to average $O(n)$ terms in
(\ref{average}) rather than $2^n$, in order to obtain a norm that is
isometric to $||| \cdot |||$, and in particular symmetrize our
original norm $\| \cdot \|$ with selected vectors $\{v_i\}_{i=1}^n$
to become almost unconditional.

In this paper, we change the settings given in (\ref{blm}), in two
ways:

\noindent First, random sign vectors are replaced by random vectors
$a(1), \ldots ,a(N)$ with an arbitrary log-concave probability
measure $\mu$.

\noindent Second, we are now interested in using a small number
$N=(1+\delta)n$ of random vectors, at the cost of an isomorphic,
instead of an almost isometric, comparison of $ |||x||| =
\int_{\R^n} \|\sum_{i=1}^n a_i x_i v_i\|d\mu(a)$ with the random
norm $ ||| x |||_N = {1\over N}\sum_{j=1}^N\|\sum_{i=1}^n a(j)_i x_i
v_i\|$.

\vskip 10pt

Let $v_1, \ldots ,v_n$ be unit vectors in a normed space $(\R^n,\|
\cdot \|)$. Define a norm $||| \cdot |||$ on $\R^n$: $||| x ||| =
\int_{\R^n} \|\sum_{i=1}^n a_i x_i v_i\|d\mu(a)$, where $a_i$ is the
$i^{th}$ coordinate of the vector $a$ and $\mu$ is a log-concave
probability measure $\mu$ on $\R^n$.

\begin{thm*}
Let $N = (1 + \delta) n$, $0<\delta<1$, and let $\{ a(1), \ldots,
a(N) \in R^n \}$ be a set of $N$ independent random vectors,
distributed with respect to a log-concave probability measure $\mu$
on $\R^n$. Then
\begin{equation*}
\P \left\{ \forall x \in \R^n \,\, c(\delta) ||| x ||| \leq
\frac{1}{N} \sum_{j=1}^N \| \sum_{i=1}^n a(j)_i x_i v_i \| \leq C
||| x ||| \right\} \ge 1-e^{-c' n} ~,
\end{equation*}
where $c(\delta)=(c\delta)^{1+{2\over\delta}}$ and $c,c',C>0$ are
universal constants.
\end{thm*}

\begin{rmk*}
It is easy to see that once you learn the theorem for small
$\delta$, it holds for large $\delta$ as well. Thus we may always
assume that $\delta<\delta_0$ for some universal constant
$\delta_0$.
\end{rmk*}

\begin{rmk*}
The norm $||| \cdot |||$ depends on the choice of vectors the $v_1,
\ldots ,v_n \in \R^n$. Different choice of vectors defines a
different norm for $x\in\R^n$. For example, if we choose $v_1=
\ldots =v_n$, it reduces the question above to the scalar case, the
norm $||| \cdot |||$ is the Euclidean norm and we get a result of
isomorphic Khinchin-type inequality; these type of results were
proved by Litvak, Pajor, Rudelson, Tomczak-Jaegermann and Vershynin
\cite{lprt},\cite{lprtv} and also by Artstein-Avidan, Friedland and
Milman \cite{afm2}.
\end{rmk*}

\begin{rmk*}
The question above was also investigated for any $0<p<\infty$. See,
e.g. Bourgain \cite{b}, Rudelson \cite{r}, Giannopoulos and Milman
\cite{gm}, Gu\'{e}don and Rudelson \cite{gr}. We shall focus on the
case of $p=1$.
\end{rmk*}

\begin{ack*} I thank my supervisor Vitali Milman for encouragement, support and interest in this work.
I thank Sasha Sodin and Shiri Artstein-Avidan for useful discussions
and remarks on preliminary versions of this note. I thank Sergey
Bobkov for suggesting me to look at the continuous case of this
problem, and last but not least, I thank Apostolos Giannopoulos for
useful remarks and referring me to the paper of Lata{\l}a \cite{la}.
\end{ack*}

\section{Proof of Theorem}
Before we proceed, we give a short description of Chernoff's method.
The following lemma, which is a version of Chernoff's bounds, gives
estimates for the probability that at least $\beta N$ trials out of
$N$ succeed, when the probability of success in one trial is $p$
(cf. Hagerup and R¨ub \cite{hr}).

\begin{lem}[Chernoff]\label{chernoff}
Let $Z_1, \ldots ,Z_N$ be independent Bernoulli random variables
with mean $0<p<1$, that is, $Z_i$ takes value $1$ with probability
$p$ and value $0$ with probability $(1-p)$. Then we have

\noindent 1)~~~~~$\P \{Z_1+ \cdots +Z_N \ge \beta N \} \ge 1-e^{-N
I(\beta,p)}$~~~~~when $\beta<p$ ~,

\noindent 2)~~~~~$\P \{ Z_1+ \cdots +Z_N \ge \beta N \} \le e^{-N
I(\beta,p)}$~~~~~~~~~~when $\beta>p$ ~,

where $I(\beta,p)=\beta\ln{\beta\over p}+(1-\beta)\ln{1-\beta\over
1-p}$.
\end{lem}

In the questions above, essentially, we are looking for upper and
lower bounds of ${1\over N} \sum_{j=1}^N\|\sum_{i=1}^n a(j)_i x_i
v_i\|$. Upper bounds are relatively easy to obtain, and quite often
do not require new methods, but only the use of large deviation
inequalities like Bernstein's inequality and some net argument.
Obtaining lower bounds is different, usually one needs small ball
probabilities, which are hard to get, and some extra delicate
arguments which are closely related to the context of the question
at hand. Here comes Chernoff's method, if one has a small ball
probability for one trial, using Chernoff's bounds, the estimate of
the average of many trials can be amplified. For more detailed
description of this method, see Artstein-Avidan, Friedland and
Milman \cite{afm1}, \cite{afm2}.

Finally, let us analyze the function
\[ I(\beta,p) = \beta\ln{\beta\over p}+(1-\beta)\ln{1-\beta\over
1-p} = u(\beta)-\beta\ln p-(1-\beta)\ln{(1-p)} ~,\] where we denoted
$u(\beta) = \beta\ln{\beta} + (1-\beta)\ln{(1-\beta)}$. The term
$u(\beta)$ is a negative, convex function which approaches $0$ as
$\beta \to 0$ and as $\beta \to 1$, and is symmetric about $1/2$
where it has a minima equal to $-\ln 2$. Thus the whole exponent in
Lemma \ref{chernoff} is of the form
\begin{equation}\label{analyse}
e^{-N I(\beta,p)}=p^{\beta
N}(1-p)^{(1-\beta)N}e^{-Nu(\beta)}\le(1-p)^{(1-\beta)N}2^N ~.
\end{equation}

\vskip 10pt \noindent \textbf{Proof of Theorem.} We estimate the
following probability
\[\P \left \{ a(1), \ldots ,a(N):\forall x\in\S_{||| \cdot |||}~~c(\delta) \le ||| x |||_N \le
C \right \} ~,\] where $||| x |||_N = {1\over
N}\sum_{j=1}^N\|\sum_{i=1}^n a(j)_i x_i v_i\|$ and $\S_{||| \cdot
|||}=\{x\in\R^n:||| x |||=1\}$.

\begin{rmk*}
The norm $||| x |||_N$ is a {\em random} norm depending on the
choice of $N$ random vectors $a(1), \ldots ,a(N)$.
\end{rmk*}

\noindent It is clear that the probability above is greater than the
following one
\begin{eqnarray*}
1&-&\P \{ \exists x\in\S_{||| \cdot |||}, ||| x |||_N>C \} \\
&-&\P \{ (\forall y\in\S_{||| \cdot |||} , ||| y |||_N\le
C)~and~(\exists x\in\S_{||| \cdot |||} , ||| x |||_N < c(\delta))
\}.
\end{eqnarray*}

\noindent \textbf{Upper bound:} We begin by estimating the first
term \[ \P \{ \exists x \in \S_{||| \cdot |||} , ||| x |||_N > C \}
~ .\] This is relatively easy, and does not require a new method; we
do it in a similar way to the one in \cite{blm}: let ${ \cal N } =
\{ y(i) \}_{ i = 1 } ^m$ be a ${1 \over 2}$-net with respect to $|||
\cdot |||$ on $ \S_{||| \cdot |||}$, it is known that such a net
exists with $m \le 5^n$. For each $1 \le i \le m$ we consider the
random variables $\{X_{i,j}\}_{j=1}^N$ defined by
\[X_{i,j}=\|\sum_{k=1}^n a(j)_k y(i)_k v_k\|-r ~,\]
where $r=\E\|\sum_{k=1}^n a_k x_k v_k\|=1$ and $y(i)_k$ is the
$k^{th}$ coordinate of the vector $y(i)$. Clearly, for each $i$ and
$j$, $X_{i,j}$ has mean $0$ and $\|X_{i,j}\|_{\psi_1} \le b$ for
some absolute constant $b>0$ (see Milman and Schechtman \cite{ms},
App. III). Now, using the well-known Bernstein's inequality which we
shall use in the form of $\psi_1$ estimate (see Vaart and Wellner
\cite{vw}):

\begin{lem}[Bernstein]\label{bernstein}
Let $Y_1, \ldots, Y_N$ be independent random variables with mean $0$
such that for some $b>0$ and every $i$, $\|Y_i\|_{\psi_1} \le b$.
Then, for any $t>0$,
\begin{equation*}
\P \{ | {1 \over N} \sum_{i=1}^N Y_i | >t \} \le 2 \exp ( -cN \min (
{t \over b},{ t^2 \over b^2}) )~,
\end{equation*}
where $c>0$ is an absolute constant.
\end{lem}

we deduce that for any $t>0$ and every $1\le i \le m$, we have
\[ \P \{ a(1), \ldots, a(N) : |{1\over N} \sum_{j=1}^N X_{i,j}|>t \} \le 2 \exp ( -cN \min (
{t \over b},{ t^2 \over b^2}) ) ~ ,\] which implies that for a point
$y(i)\in\cal{N}$ and for any $t>1$ we have
\begin{eqnarray*}
\P \big \{ a(1), \ldots , a(N) : {1 \over N} \sum_{j=1}^N \|
\sum_{k=1}^n a(j)_k y(i)_k v_k \| > t \big \} \le \\
\le 2 \exp ( -c N \min ( {t-1 \over b}, { (t-1)^2 \over b^2} )).
\end{eqnarray*}

The obvious way to make this probability small enough to handle a
large net is to increase $t$, and obviously we shall get a worse
upper bound constant. So, we choose $t$ such that
\[2 e^ { -c N \min ( {t-1 \over b}, { (t-1)^2 \over b^2} ) } \cdot 5^n\le e^{-n} ~,\] for example
$ t = 6 \max({b \over \sqrt c}, {b \over c}) + 1 $. Then, with
probability at least $1-e^{-n}$, for every $1 \le i \le m$ we have
\[{1\over N} \sum_{j=1}^N \|\sum_{k=1}^n a(j)_k y(i)_k v_k\| \le t ~.\] We thus have an upper bound for a net on the sphere. It is
standard to transform this to an upper bound estimate on {\em all}
the sphere (this is an important difference between lower and upper
bounds). One uses a consecutive approximation of a point on the
sphere by points from the net to get that $||| x |||_N\le 2t = 2t$
for every $x \in \S_{||| \cdot |||}$. This completes the proof of
the upper bound, where $ C = 12\max({b \over \sqrt c}, {b \over c})
+ 2$ is a universal constant.

\vskip 10pt \noindent \textbf{Lower bound:} We now turn to estimate
the second term
\begin{equation}\label{lower}
\P \{ (\forall y\in\S_{||| \cdot |||} ,|||y|||_N\le C)~and~(\exists
x\in\S_{||| \cdot |||} ,|||x|||_N<c(\delta)) \} .
\end{equation}

Note that when estimating this term, we know {\em in advance} that
the (random) norm $|||\cdot|||_N$ is bounded from above on the
sphere $\S_{|||\cdot|||}$ (i.e. $\forall y\in\S_{||| \cdot |||}$ we
have $|||y|||_N\le C$, where $C$ comes from the upper bound). This
is crucial to transform a lower bound on a net on the sphere to a
lower bound on the whole sphere. For the lower bound we use
Chernoff's method, as described above, to estimate the probability
in (\ref{lower}).

Let us denote by $p$ the probability that for a random vector
$a\in\R^n$ we have $\|\sum_{i=1}^n a_i x_i v_i\| \ge \alpha$, where
$\alpha>0$ and $x$ is some point on $\S_{|||\cdot|||}$:
\[p=\mu(\{a\in\R^n:\|\sum_{i=1}^n a_i x_i v_i\| \ge \alpha\}) ~.\]

If "doing an experiment" means checking whether $\|\sum_{i=1}^n a_i
x_i v_i\| \ge \alpha$ (where $a\in\R^n$ is a random vector) then for
$|||x|||_N$ to be greater than some $c$, it is enough that
$(c/\alpha)N$ of the experiments succeed.

Of course, we will eventually not want to do this on \textit{all}
points $x$ on the sphere, but just on some dense enough set. So,
first we estimate the probability $p~$:

\begin{lem}\label{simple}
There exists a universal constant $\gamma>0$ such that for any
$x\in\S_{||| \cdot |||}$, we have
\[\mu(\{a\in\R^n:\|\sum_{i=1}^n a_i x_i v_i\|\le\gamma\})\le{2\over3} ~.\]
\end{lem}

\noindent \textbf{Proof of Lemma \ref{simple}.} Let us define
$A_x=\{a\in\R^n:\|\sum_{i=1}^n a_i x_i v_i\|\le\gamma_x\}$ for
$x\in\S_{||| \cdot |||}$, then $A_x$ is convex and symmetric set. We
take $\gamma_x>0$ to be the number such that $\mu(A_x)={2 \over 3}$.
Applying Borell's lemma (see Milman and Schechtman \cite{ms}, App.
III.3) we get, for all $t>1$,
\[\mu (\{a\in\R^n:\|\sum_{i=1}^n a_i x_i v_i\|>\gamma_x t\})\le {2\over 3}({1\over
2})^{(t+1)/2} ~,\] and consequently,
\begin{eqnarray*}
1=|||x|||&=&\int_{R^n}\|\sum_{i=1}^n a_i x_i v_i\|d\mu(a)=
\int_0^\infty \mu(\|\sum_{i=1}^n a_i x_i v_i\|>t)dt\\
&\le& \gamma_x + \int_{\gamma_x}^\infty \mu(\|\sum_{i=1}^n a_i x_i
v_i\|>t)dt\\
&\le& \gamma_x + \gamma_x \int_1^\infty
{2\over3}({1\over2})^{t+1\over 2}dt\\
&\le& \gamma_x (1+c_4)
\end{eqnarray*}

for some universal constant $c_4>0$ which doesn't depend on $x$.
Therefore, we get that $\gamma_x\ge{1\over{1+c_4}}>0$. Now, we take
$\gamma={1\over{1+c_4}}$. For this $\gamma>0$ and any $x\in\S_{|||
\cdot |||}$ we have $\mu(\{a\in\R^n:\|\sum_{i=1}^n a_i x_i
v_i\|\le\gamma\})\le {2\over 3}$ $\hfill\square$

\vskip 10pt \noindent Now, we can use the following lemma of
Lata{\l}a with the set \[C_x=\{a\in\R^n:\|\sum_{i=1}^n a_i x_i
v_i\|\le \gamma\} ~,\] for which we proved above that $\mu(C_x)\le{2
\over 3}=b<1$.

\begin{lem}[Lata{\l}a \cite{la}]\label{latala}
For each $b<1$ there exists a constant $c_b>0$ such that for every
log-concave probability measure $\mu$ and every measurable convex,
symmetric set $C$ with $\mu(C)\le b$ we have \[\mu(tC)\le c_b t
\mu(C)~~~~~for~t\in[0,1] ~.\]
\end{lem}

Notice that, for \textit{any} point $x$, we can make the probability
as small as we like by reducing $t$. This allows us to use a simple
net: take $\theta$-net ${\cal N}$ in $\S_{|||\cdot|||}$, with less
than $({3\over \theta})^n$ points.

For every $x\in\S_{|||\cdot|||}$ there is a vector $y\in{\cal N}$
such that $|||x-y|||\le\theta$, and we have
$|||y|||_N\le|||x|||_N+|||x-y|||_N\le c+C\theta$ (where
$c=c(\delta)$ and $C$ comes from the upper bound). Therefore, we
bound (\ref{lower}) by
\begin{equation}\label{cont_net}
\P \{ \exists y\in{\cal N}: |||y|||_N\le c+C\theta \}.
\end{equation}

By Lemma \ref{latala}, for a given $y\in{\cal N}$ we have for any
$0<t<1$ that
\[p=\mu(\{a\in\R^n:\|\sum_{i=1}^n a_i y_i v_i\|\ge t \gamma\})\ge
1-c_1 t ~,\] where $c_1=c_b{2\over3}$.

We return to our scheme, in order to estimate the probability in
(\ref{cont_net}), assume that $\beta t \gamma \ge c+C\theta$, where
$\beta$ shall be the portion of good trials out of $N$, and $t$
another constant that we choose later such that $p>\beta$. So, we
know that for $\beta<1-c_1 t$ (which is hardly a restriction, $t$
will be very small and so will $\beta$), from Lemma \ref{chernoff}
for a given $y\in{\cal N}$ we have
\begin{eqnarray*}
\P \{ |||y|||_N \ge \beta t \gamma \} \ge \P \{Z_1+\cdots+Z_N \ge
\beta N \} \ge 1-e^{-NI(\beta,p)} ~,
\end{eqnarray*}
where $\P\{Z_i=1\}=\mu(\{a\in\R^n:\|\sum_{i=1}^n a_i y_i v_i\|\ge t
\gamma\})$.

We choose $\beta$ so that $(1+\delta)(1-\beta)=1+{\delta\over 2}$,
hence $\beta={\delta\over {2(1+\delta)}}$. We choose $\theta=\beta t
\gamma /2C$, where $C$ comes from the upper bound. To make sure that
the probability above holds for all points in the net we ask that
\[ e^{-NI(\beta,p)}\cdot|{\cal{N}}|\le 2^N
(1-p)^{(1-\beta)N}\cdot({3\over\theta})^n \le \]
\[ \le 2^N (c_1 t)^{(1+{\delta\over2})n}\cdot({3\over\theta})^n \le
{1\over 2}e^{-c'n} ~.\] For the first inequality we use
(\ref{analyse}), we choose $t=(c_2\delta)^{2\over\delta}$, for some
universal constant $c_2>0$, and get the lower bound for each point
of the net ${\cal N}$.

Now using the upper bound, for every $x\in\S_{|||\cdot|||}$ there is
a vector $y\in\cal{N}$ such that $|||x-y|||\le\theta$, therefore we
have \[|||x|||_N \ge ||| y |||_N - ||| x-y |||_N \ge \beta t
\gamma-\theta C=:c ~,\]
$c=c(\delta)=(c_3\delta)^{1+{2\over\delta}}$, where $c_3>0$ is an
absolute constant. Thus the proof of the lower bound, and of the
Theorem, is completed. $\hfill\square$


\begin{thebibliography} {GGM}

\bibitem {afm1} S. Artstein-Avidan, O. Friedland, V.D. Milman. {\em
Geometric applications of Chernoff-type estimates and a zigzag
approximation for balls}. Proc. Amer. Math. Soc. 134 (2006), no. 6,
1735--1742.

\bibitem {afm2} S. Artstein-Avidan, O. Friedland, V.D. Milman. {\em
Geometric Applications of Chernoff-Type Estimates}, Geometric
Aspects of Functional Analysis, Israel Seminar 2004--2005, Lecture
Notes in Mathematics , Vol. 1910.

\bibitem {b} J. Bourgain. {\em Random points in isotropic
convex sets}. Convex Geometric Analysis (Berkeley, CA, 1996), 53-58,
Math. Sci. Res. Inst. Publ., 34, Cambridge Univ. Press, Cambridge,
1999.

\bibitem {blm} J. Bourgain, J. Lindenstrauss, V.D. Milman.
{\em Minkowski sums and symmetrizations}. Geometric aspects of
functional analysis (1986/87), Lecture Notes in Math., 1317,
Springer, Berlin, 1988, 44-66.

\bibitem {gm} A.A. Giannopoulos, V.D. Milman. {\em
Concentration property on probability spaces}. Adv. Math. 156
(2000), no. 1, 77--106.

\bibitem {gr} O. Gu\'{e}don, M. Rudelson. {\em $L\sb
p$-moments of random vectors via majorizing measures}. Adv. Math.
208 (2007), no. 2, 798--823.

\bibitem {hr} T. Hagerup, C. R¨ub. {\em A guided tour of
Chernoff bounds}. Info Proc Lett 33 (1990) no. 6, 305-308.

\bibitem {ka} J. P. Kahane. {\em Some random series of functions}. second edition, Cambridge Studies in Advanced Math., 5, (1985),
Cambridge University press.

\bibitem {kw} S. Kwapie\'n. {\em A theorem on the Rademacher
series with vector valued coefficients}. Probability in Banach
spaces, Oberwolfach 1975, Lecture Notes in Mathematics, Vol. 526,
Springer-Verlag, Berlin-Heidelberg-New York, 1976, pp. 157–158.

\bibitem {la} R. Lata{\l}a. {\em On the equivalence between
geometric and arithmetic means for log-concave measures}. in Convex
geometric analysis (Berkeley, CA, 1996), 123-127, Math. Sci. Res.
Inst. Publ. 34, Cambridge Univ. Press, Cambridge 1999.

\bibitem {lprt} A.E. Litvak, A. Pajor, M. Rudelson, N. Tomczak-Jaegermann. {\em Smallest singular value of random
matrices and geometry of random polytopes}. Adv. Math. 195 (2005),
no. 2, 491--523.

\bibitem {lprtv} A.E. Litvak, A. Pajor, M. Rudelson, N.
Tomczak-Jaegermann, R. Vershynin. {\em Random Euclidean embeddings
in spaces of bounded volume ratio}. C. R. Math. Acad. Sci. Paris 339
(2004),  no. 1, 33--38.

\bibitem {ms} V.D. Milman and G. Schechtman. {\em Asymptotic
theory of finite dimensional normed spaces}. Springer Lecture Notes,
1200, (1986).

\bibitem {r} M. Rudelson. {\em Random vectors in the isotropic
position}. J. Funct. Anal. 164 (1999), no. 1, 60--72.

\bibitem {vw} A.W.v.d. Vaart, J.A. Wellner. {\em Weak Convergence and Empirical Processes}. Springer Series
in Statistics, 2000.

\end{thebibliography}
\end{document}